\documentclass[opre]{informs3}

\OneAndAHalfSpacedXI % current default line spacing
\usepackage{natbib}
 \bibpunct[, ]{(}{)}{,}{a}{}{,}%
 \def\bibfont{\small}%
 \def\bibsep{\smallskipamount}%
\DeclareSymbolFont{matha}{OML}{txmi}{m}{it}% txfonts
\DeclareMathSymbol{\varv}{\mathord}{matha}{118}

\usepackage{bbm}
\def\N{\mathbb{N}}
\def\R{\mathbb{R}}

\def \E {\mathbb{E}}
\def \P {\mathbb{P}}

\def \C {\mathcal{C}}
\def \-> {\rightarrow}

\def \xb {\mathbf{x}}

\def \eqd {\,{\buildrel d \over =}\,}

\def \Var {\mbox{\bf Var}}

\def \Jinf {J_{\infty}}

\def \tT {\Tilde{t}}
\def \PC {\mathcal{P}}
\def \xiT {\Tilde{\xi}}

\usepackage{enumerate}
\usepackage[shortlabels]{enumitem}
\usepackage{tabularx,booktabs,mathrsfs,setspace}%
\usepackage{mathtools}
\usepackage{multirow}
\usepackage{rotating}
\usepackage{float}
\usepackage{pgfplots}
\pgfplotsset{compat=1.9}
\usepackage{csvsimple}
\usepackage{makecell}
\usepackage{multirow}
\usepackage{adjustbox}
\usepackage[caption=false, position=top]{subfig}
\usepackage{alphalph}

\usepackage{mathtools}
\usepackage{comment}
\usepackage{lscape}
\usepackage{longtable}

\TheoremsNumberedThrough     % Preferred (Theorem 1, Lemma 1, Theorem 2)
\EquationsNumberedThrough    % Default: (1), (2), ...
\begin{document}
%%%%%%%%%%%%%%%%

%\setcounter{page}{1}
%\DOI{10.1287/xxxx.0000.0000}%

\RUNAUTHOR{Moradi et al.} % for four or more authors

\RUNTITLE{Asymptotic Optimality of Projected Inventory Level Policies}

% \TITLE{Bundling Information Goods of Decreasing Value}
\TITLE{Asymptotic Optimality of Projected Inventory Level Policies for Lost Sales Inventory Systems with Large Leadtime and Penalty Cost}

\ARTICLEAUTHORS{%
\AUTHOR{Poulad Moradi, Joachim Arts}
\AFF{Luxembourg Centre for Logistics and Supply Chain Management, University of Luxembourg, Luxembourg-City, Luxembourg, \\ \EMAIL{\{poulad.moradi, joachim.arts\}@uni.lu}}
\AUTHOR{Melvin Drent}
\AFF{Department of Information Systems and Operations Management, Tilburg University, the Netherlands, \\ \EMAIL{m.drent@tilburguniversity.edu}
}}% end of the block

\ABSTRACT{%
We study the canonical periodic review lost sales inventory system with positive leadtime and independent and identically distributed (i.i.d.) demand under the average cost criterion. We demonstrate that the relative value function under the constant order policy satisfies the Wiener-Hopf equation.
We employ ladder processes associated with a random walk featuring i.i.d. increments, to obtain an explicit solution for the relative value function. 
This solution can be expressed as a quadratic form and a term that grows sublinearly. Then we perform an approximate policy iteration step on the constant order policy and uniformly bound the gap relative to the otimal cost rate for large lead times. This leads to our main result that projected inventory level policies are asymptotically optimal as the leadtime grows when the cost of losing a sale is sufficiently large and demand has a finite second moment.
% Under these conditions, we also show that the optimal cost rate approaches infinity, proportional to the square root of the cost of losing a sale.
}%

%\KEYWORDS{deterministic inventory theory; infinite linear programming duality; 
%  existence of optimal policies; semi-Markov decision process; cyclic schedule}

\KEYWORDS{Lost Sales, Asymptotic Optimality, Markov Decision Processes, Inventory}

\maketitle
%%%%%%%%%%%%%%%%%%%%%%%%%%%%%%%%%%%%%%%%%%%%%%%%%%%%%%%%%%%%%%%%%%%%%%

% Samples of sectioning (and labeling) in TRSC
% NOTE: (1) \section and \subsection do NOT end with a period
%       (2) \subsubsection and lower need end punctuation
%       (3) capitalization is as shown (title style).
%
%\section{Introduction.}\label{intro} %%1.
%\subsection{Duality and the Classical EOQ Problem.}\label{class-EOQ} %% 1.1.
%\subsection{Outline.}\label{outline1} %% 1.2.
%\subsubsection{Cyclic Schedules for the General Deterministic SMDP.}
%  \label{cyclic-schedules} %% 1.2.1
%\section{Problem Description.}\label{problemdescription} %% 2.

% Text of your paper here
\section{Introduction}
The control of lost sales inventory systems remains a fundamental challenge in inventory theory.
In such systems, unmet demand caused by stockouts is lost, often resulting in substantial penalty costs.
We consider the canonical lost sales inventory system, which is a single-item, single-echelon, periodic-review inventory system with a positive leadtime and independent and identically distributed (i.i.d.) demand under the average cost criterion.
This system serves as the foundation for more complex lost sales inventory models.
Therefore, developing well-performing and computationally efficient control policies for the canonical system is crucial to derive effective policies for real-world lost sales inventory problems.

The optimal replenishment policy for the canonical system with negligible leadtime reduces to a newsvendor problem. When the leadtime is positive, the optimal policy can be found through dynamic programming but this is intractable due to the curse of dimensionality.
Consequently, a key inventory research stream in stochastic lost sales inventory control focuses on developing simple heuristic policies that perform well under specific conditions, such as achieving asymptotic optimality in certain scaling regimes.
We refer interested readers to \cite{Bijvank2023}, \cite{Goldberg2021}, and \cite{Bijvank2011}, for further discussions on lost sales inventory systems and related asymptotic optimality results.

\cite{Huh2009} and \cite{Bijvank2014} analyze base-stock policies that place orders to raise the inventory position to a fixed base-stock level.
They establish that such policies are asymptotically optimal as the cost of losing a sale grows for a fixed leadtime.
\cite{Goldberg2016} and \cite{Xin2016} demonstrate that a constant order policy, which places the same order quantity every period, is asymptotically optimal as the leadtime grows for a fixed cost of losing a sale.
Both the base-stock policy and the constant order policy rely on a single parameter, making them easy to implement in practice. However, neither policy is optimal across both asymptotic regimes.
To address this limitation, \cite{Xin2021} proposes a two-parameter hybrid policy that integrates the base-stock and constant order policies, and proves its asymptotic optimality for large leadtimes.
This policy, known as the capped base-stock policy, was initially studied by \cite{Johansen2008} and can be readily shown to be asymptotically optimal as the cost of losing a sale grows large.
By adjusting its parameters, the capped base-stock policy can thus be tailored to achieve asymptotic optimality in either regime.

Recently, \cite{vanJaarsveld2024} introduced the projected inventory level (PIL) policy, which places orders to ensure that the expected inventory level at the time of receipt reaches a fixed target.
Unlike constant order and base-stock policies, the PIL policy dynamically adjusts order quantities by leveraging probabilistic information available at each decision epoch.
\cite{vanJaarsveld2024} demonstrate that the PIL policy consistently outperforms the base-stock policy for general demand distributions and prove that it also outperforms the constant order policy when demand is exponential. PIL policies are also asymptotically optimal for perishable inventory systems in several regimes \citep{Bu2025a,Bu2025b}, and the projection idea is similarly employed by \citet{drent2022effective} for dual-sourcing inventory systems, where it yields both asymptotic optimality and strong empirical performance.

Policies developed for the canonical lost sales system can be extended to more complex settings, including systems with non-stationary demand, perishable items, continuous review, partially observable parameters, finite storage capacity, supply uncertainty, stochastic returns, joint inventory and pricing control, and finite horizon decision making  \citep[see, e.g.,][]{Bu2025a, Bu2025b, Bu2024,  Bu2020, Lyu2024, Bai2023, Xin2022, Chen2021}.

The asymptotic regimes discussed in the literature only consider one parameter growing large while keeping all other parameters fixed. 
In this paper we study the performance of the PIL policy under a general demand process as the leadtime grows large when the cost of losing a sale is sufficiently large.
Under mild conditions on the demand distribution we show that:
\begin{enumerate}
    \item The  difference between the average cost-rate of the PIL policy and that of the constant order policy remains bounded by a finite constant .
    \item The PIL policy is asymptotically optimal for sufficiently large lost sales penalty costs as the leadtime approaches infinity.
\end{enumerate}
Our analysis hinges on new bounds we derive for the solution of the Wiener-Hopf equation that characterizes the relative value function under the constant order policy. These bounds follow from studying the ladder processes of a random walk with increments equal to the per-period excess demand minus the constant order. We then apply a one-step policy improvement technique to analyze the cost-rate difference between the PIL and constant order policies.

The rest of the paper is organized as follows. Section \ref{section:Model_MainResult} introduces the model and optimization problem (Section \ref{section:Model}), as well as the main result (Section \ref{section:MainResult}).
Section \ref{section:MainResultProof} provides the proof of the main result, including the introduction of ladder processes (Section \ref{section:LadderProcesses}), the solution to our Wiener-Hopf equation (Section \ref{section:WHequation}), asymptotic inventory dynamics (Section \ref{section:InventoryDynamicsAsymp}), and policy improvement argument (Section \ref{section:additionalCostRate}) which completes the proof.
A summary of results and final remarks are provided in Section \ref{Section:Conclusion}.
All proofs are included in the Appendix, unless otherwise specified.

\section{Model and main result}
\label{section:Model_MainResult}
\subsection{Model}
\label{section:Model}

% \begin{table*}[h!]
%   \caption{\textsf{Overview of notation.}}
% 	\fontsize{8pt}{9pt}\selectfont
% \scriptsize
%   \label{tab:notation}
%   \begin{tabularx}{\textwidth}{l X}
%     \toprule
%     Notation& Description\\
% 		\midrule
% 	Sets & \\
%         \\ \\
%     Parameters &  \\
%         $h$             & Unit holding cost.\\
%         $p$             & Unit lost sales penalty cost.\\
%         $L$             & Deterministic leadtime.\\
%         $D_t$           & Demand at time $t$, i.i.d., distribute as the non-negative r.v. $D \sim \mathcal{D}$.\\
    
% 	    \\ \\
% 	Variables & \\
%         $q_t$           & Order quantity received at time $t$.\\
%         $J_t^{\pi}$     & Inventory level at the beginning of time $t$ before the order arrival under the policy $\pi$.\\
%         \\
%     \bottomrule
%   \end{tabularx}
% \end{table*}
We consider an infinite-horizon periodic review lost sales inventory system. Demand in period $t$ is denoted $D_t$ and $\{D_t\}_{t\in \N_0}$ ($\N_0:=\N \cup \{0\}$) is a sequence of non-negative independent and identically distributed random variables with distribution function $F_D$ supported on $[0,\infty)$, and finite mean $\mu_D := \E[D] <\infty$ and variance $\Var[D] :=\sigma_D^2 \in (0,\infty)$. We assume $F_D(0) = 0$ for notational simplicity, though all results remain valid without this assumption.
Each time period $t \in \N_0$ we receive an order, $q_t \in \R_+$, that is placed in period $t-L$, where $L \in \N_0$ is the deterministic leadtime.
Let $\{J_t\}_{t \in \N_0}$ denote the sequence of inventory level random variables at the beginning of each period before receiving the order.
The state of the system at the start of period $t \in \N_0$, denoted by $\xb_t \in \R_+^{L+1}$, is a vector comprising the inventory level in period $t$ as well as the outstanding orders in the pipeline.
That is, $\xb_t = (J_t, q_t, q_{t+1}, \dots, q_{t+L-1})$.
We assume that $\xb_0$ is fixed and known, and $J_0=0$.
Demand that exceeds the on-hand inventory $J_t+q_t$ is lost at the end of the period at the unit cost of $p\geq0$. Any surplus inventory at the end of a period is held at a cost of $h\geq0$ per item. 
The sequence of events in each period $t \in \N_0$ is as follows: 
(1) The state of the system $\xb_t$ is observed and the order $q_{t+L}$ is placed, (2) The order $q_t$ is received, (3) The demand $D_t$ is realized, and (4) the costs of period $t$ are incurred as $p(D_t-q_t-J_t)^+ + h (J_t+q_t-D_t)^+$ where where, $(x)^+ := \max(x, 0)$. The dynamics of the inventory level are
\begin{equation}
    \label{eq:JDynamics}
    J_{t+1} = (J_t + q_t - D_t)^+.
\end{equation}

A policy $\pi$ is a set of mappings from the space of the states, $\xb_t$, to the space of orders, $q_{t+L}$, i.e., $\{\pi_t: \R_+^{L+1} \to\R_+\}_{t \in \N_0}$. We denote by $\Pi$ the set of admissible policies.
A policy $\pi$ is stationary if $\pi_t(\xb)=\pi_0(\xb)$ for all $t \in \N_0$ and $\xb \in \R^{L+1}$.
When a policy $\pi$ is stationary, we omit the index $t$ in $\pi_t$, for simplicity.
We denote by $q_t(\pi)$ and $J_t(\pi)$ the random variables for the order quantity and inventory level respectively under policy $\pi$.
We consider two stationary policies: the constant order policy $C_r$ \citep[cf.][]{Xin2016} and the projected inventory level (PIL) policy $P_{\xi}$ \citep[cf.][]{vanJaarsveld2024}, where $r \in [0,\mu_D)$ is the constant order and $\xi\geq0$ is the projected inventory level.
For $t \in \N_0$, the constant order policy and PIL policy are expressed by:
\begin{equation*}
    C_r(\xb_t):=r, \quad \mbox{ and } \quad
    P_{\xi}(\xb_t):= (\xi - \E[J_{t+L}\vert \xb_t])^+.
\end{equation*}
Let $\{\{c_t(\pi)\}_{t\in \N_0}\}_{\pi \in \Pi}$ be the sequence of cost random variables given by:
\[
    c_t (\pi) := h \left(J_t(\pi) + q_t(\pi) - D_t \right)^+ + p \left(D_t - J_t (\pi) - q_t(\pi) \right)^+.
\]
As a notational convenience, we define $D_{[a,b]} = \sum_{t=a}^b D_t$, and similarly define $J_{[a,b]}$, $q_{[a,b]}$, and $c_{[a,b]}(\pi)$.
Accordingly, the cost-rate function $\C: \Pi \to \R_+$ is defined as:
\[
    \C(\pi) :=
    \limsup_{T \to \infty} \E\left[\frac{c_{[L,T]}(\pi)}{T-L+1}\right].
\]
We will sometimes write the dependence of $\C(\pi)$ on $p$ and $L$ explicitly as $\C(\pi \mid p,L)$.
Let $\C^*(p,L):= \inf_{\pi \in \Pi} \C(\pi\mid p,L)$ denote the optimal cost-rate.
\cite{Huh2011} show that a stationary policy $\pi^*$ exists such that $\C(\pi^*) = \C^*$.
Throughout the paper we say a function $g$ is $o(f(x))$ and write $g(x)=o(f(x))$ if and only if $\lim_{x\to\infty} g(x)/f(x)=0$. 

\subsection{Main result}
\label{section:MainResult}
In this section we present the main result.
For a fixed demand distribution $F_D$ and $h$, we construct a sequence $\{\xi_p\}_{p\geq0}$ such that $\xi_p \in \arg \min_{\xi\geq0} \lim_{L \to \infty}\C(P_{\xi} \mid p,L)$, and a sequence $\{r_p\}_{p\geq0}$ such that  $r_p \in \arg \min_{r\in [0,\mu_D)} \lim_{L \to \infty}\C(C_r \mid p,L)$. We now state our main result.
\begin{theorem}
    \label{theorem:epsilonOptimality}
    \begin{enumerate}
        \item There exists a constant $0 \leq M < \infty$ such that the optimality gap of the best PIL policy remains bounded by $M$ for all $p\geq 0$ as $L$ tends to infinity, that is
        \[
            \C\left(P_{\xi_p}\right) - \C(C_{r_p})=\lim_{L \to \infty} \left( \C\left(P_{\xi_p}\right) - \C^*(p,L) \right) < M \quad \mbox{for all $p \geq 0$}.
        \]
        \item The PIL policy is asymptotically optimal for large $p$ when $L$ approaches infinity, that is,
        \[
            \lim_{p \to \infty} \lim_{L \to \infty} \frac{\C\left(P_{\xi_p}\right)}{\C^*(p,L)} = 1.
        \]
    \end{enumerate}

\end{theorem}
Theorem \ref{theorem:epsilonOptimality} states that the cost-rate of the PIL policy exceeds that of the constant-order policy by at most a fixed constant $M$. Since the constant order policy is asymptotically optimal as the leadtime increases \citep{Goldberg2016,Xin2016}, the optimality gap of PIL policy remains bounded $M$ under the same condition.
% This result provides theoretical support for the strong empirical performance of the PIL policy observed in \cite{vanJaarsveld2024}, and
This result extends our understanding of the asymptotic optimality of PIL policies beyond the special case of exponentially distributed demand addressed in Theorem 2 of \cite{vanJaarsveld2024}.

We note that Theorem \ref{theorem:epsilonOptimality} differs from the asymptotic optimality result for sufficiently large $p$ presented in \cite{vanJaarsveld2024} (Theorem 4) in terms of the asymptotic regime.
The analysis in \cite{vanJaarsveld2024} relies on the comparison of the PIL policy and the base-stock policy, which is not optimal as the leadtime approaches infinity.
It is worth noting that, unlike the base-stock policy, whose optimality gap grows unboundedly with the lead time, the optimality gap of the PIL policy remains uniformly bounded.
% Second, \cite{vanJaarsveld2024} provide their optimality result in a multiplicative sense, meaning that the absolute optimality gap may not vanish in the limit as $p \to \infty$. Note that the optimal cost-rate diverges as $p$ tends to infinity. In contrast, our result shows that the absolute optimality gap of the best PIL policy approaches zero while the optimal cost-rate grows large in the asymptotic regime, where $L \to \infty$ before $p \to \infty$.

%This result raises an interesting open question 
%regarding the absolute optimality gap of the PIL policy in a regime with fixed $L$ and large $p$.
%In Section \ref{section:additionalCostRate} we will explain how $\xi_p$ relates to the best constant order policy under any $p\geq0$.

\section{Proof of Theorem \ref{theorem:epsilonOptimality}}
\label{section:MainResultProof}
Assuming exponentially distributed demand, \cite{vanJaarsveld2024} show that the relative value function of the constant-order policy is a parabola, and that a one-step policy improvement yields the PIL policy. This establishes the asymptotic optimality of the PIL policy as the leadtime grows, since it strictly improves upon the constant-order policy, which is itself asymptotically optimal \citep{Goldberg2016}. In this paper, we extend this approach to general demand distributions with finite second moments by showing that the relative value function has a quadratic form and a term that grows sublinearly. Next we use that as the cost of losing a sale grows, the optimal constant order policy approaches a heavy traffic regime where this sublinear term turns out to be unimportant and a PIL policy will not be worse than a constant order policy within tight bounds.
%However, their theoretical validation of this superiority is limited to the case of exponentially distributed demand.
%The crux of our argument hinges on a single-step improvement idea, a technique also employed by \cite{vanJaarsveld2024} for the specific case of exponential demand.
%The intuition behind this idea can be explained as follows.
%For each $p\geq0$, consider a policy $\Tilde{\pi}_p$ in which the system follows the best PIL policy $P_{\xi_p}$, i.e., order quantities $q_t(P_{\xi_p})$ are received from $L$ to $\tT\geq L$, and then transitions to the best constant order policy $C_{r_p}$, i.e., order quantities $q_t(C_{r_p})$ are received between $\tT$ and $T \geq \tT$.
%Define the long-run additional expected cost of policy $\Tilde{\pi}_p$ compared to $C_{r_p}$:
%\[
%    f(p) :=  \lim_{T \to \infty}\E[c_{[L,T]}(\Tilde{\pi}_p)-c_{[L,T]}(C_{r_p}) \vert \xb_0].
%\]
%We investigate whether for any $\epsilon>0$, there exists a threshold $p_{\epsilon}\geq0$ such that for all $p \geq p_{\epsilon}$, $f(p) \leq \epsilon$.
%While an explicit expression for $f(p)$ is not readily available, we use relative value functions under the constant order policy to show that $f(p)$ is upper bounded by a negligible value for sufficiently large $p$.
% The rigorous representation of the relative value function for the constant order policies involves a random walk argument, which we present in the next section.

Let $Y$ be a random variable with distribution $F_Y$, defined as the difference between $D$ and the constant order $r$, i.e. $Y := D - r$.
It is straightforward to verify that $F_Y$ is concentrated on $[-r, \infty)$, since $F_Y(x) = F_D (x+r)$, and $\Var[Y]=\Var[D]<\infty$. Let $\mu_Y:=\E[Y]$ and $\sigma_Y^2:=\Var[Y]$. 
\begin{definition}
    \label{def: valuefunc}
    For a constant order policy $C_r$, $r \in [0 , \mu_D)$,
    the relative value function $v_r: \R_+ \rightarrow \R$ satisfies,
    \begin{equation*}
    \label{eq:RelValinitial}
        v_r (x) := \E_Y \left[h (x-Y)^+ + p(Y-x)^+ + v_r \left((x-Y)^+\right) \right]- \C(C_r), \qquad v_r(0) = 0, \ \ x \geq 0.
    \end{equation*}
    %with the boundary condition $v_r (0)=0.$
\end{definition}

\noindent
The difference $v_r (x_1) - v_r (x_2)$ represents the additional total long-run expected cost when the system starts from $x_1$ rather than $x_2$ under the constant order policy $C_r$ \citep[cf., Chapter 6][]{Tijm2003}.
\cite{Goldberg2016} show that $\C(C_r) = h \E [\Jinf] + p \mu_Y$, where $\Jinf$ denotes the steady state inventory level under $C_r$, i.e., $\P\{\Jinf \leq x\} = \lim_{t \rightarrow \infty} \P\{J_t(C_r) \leq x\}$.
By Definition \ref{def: valuefunc}, the relative value function $v_r (x)$ can be expressed as a convolution equation:

\begin{equation}
    \label{eq:RelValint}
    v_r(x) = a_r(x) + \int_{-r}^x v_r(x-y)F_Y(dy), \qquad v_r(0) = 0, \ \ x \geq 0,
\end{equation}
where,
\begin{equation}
    \label{eq:ar}
    a_r(x) := h \E_Y[(x-Y)^+] + p \E_Y[(Y-x)^+]- p \mu_Y - h \E [\Jinf].
\end{equation}

\noindent
To simplify notation we introduce the convolution operator, $*$, as follows. Let $K: \R \rightarrow \R$ and $F: \R \rightarrow \R$ be two real functions. The convolution of $K$ and $F$ is defined as:
\[
    K*F (x) := \int_{-\infty}^x K(x-y)F(dy).
\]
Therefore, Equation \eqref{eq:RelValint} can be rewritten as,
\begin{equation}
    \label{eq:RelVal}
    v_r(x) = a_r(x) + v_r * F_Y (x), \qquad v_r(0) = 0, \ \ x \geq 0.
\end{equation}
Deriving an explicit expression for $v_r(x)$ is non-trivial, as Equation \eqref{eq:RelVal} constitutes a Wiener-Hopf equation \citep[cf.][]{Asmussen1998}. However, by analyzing a specific random walk with i.i.d. increments and its associated ladder processes in Section \ref{section:LadderProcesses}, we are able to derive an explicit solution in Section \ref{section:WHequation}, which enables the remainder of our analysis.

\subsection{Ladder processes}
\label{section:LadderProcesses}

Consider a random walk $\{S_t := \sum_{i=1}^t Y_i\}_{ t \in \N_0}$, with $S_0=0$, where $\{Y_t\}_{ t \in \N_0}$ represents a sequence of random variables defined by $Y_t := D_t -r$.
We introduce two stopping periods associated with the random walk $S_t$. The (weak) ascending ladder period, denoted by $\tau_+$, is the first period (greater than zero) that the random walk attains a non-negative value, i.e., $\tau_+ := \inf \{t >0: S_t \geq 0\}$.
The value of the stopped random walk at $\tau_+$, i.e., $S_{\tau_+}$, is a random variable known as the first (weak) ascending ladder height, with distribution function $G_+ (x) = \P \{S_{\tau_+} \leq x \}$ supported on $[0, \infty)$.
The mean and variance of $S_{\tau_+}$ are denoted by $\mu_+ := \E[S_{\tau_+}]$ and $\sigma_+^2 := \Var[S_{\tau_+}]$, respectively.
Both $\mu_+$ and $\sigma_+$ are finite for $r \in [0, \mu_D)$, and remain so as $r$ approaches $\mu_D$.
% Define the sequences of $\{\mu_+(r)\}_{r \in [0,\mu_D)}$, and $\{\sigma_+(r)\}_{r \in [0,\mu_D)}$ the mean and standard deviation of ascending ladder hight under increments $\{Y_t()\}$.

\begin{lemma}
    \label{lemma:momentsPlus}
    $\lim_{r \uparrow \mu_D} \mu_+$ and $\lim_{r \uparrow \mu_D} \sigma_+$ exist, and (a) $0 < \lim_{r \uparrow \mu_D} \mu_+ < \infty$, and (b) $\lim_{r \uparrow \mu_D} \sigma_+ < \infty$.
\end{lemma}

Likewise, the (strict) descending ladder period is the first period (greater than zero) that the random walk takes a negative value, i.e., $\tau_- := \inf \{t >0: S_t < 0\}$. Accordingly, $S_{\tau_-}$, is the first (strict) descending ladder height random variable with  distribution function $G_- (x) := \P \{S_{\tau_-} \leq x \}$ supported on $(-\infty, 0)$.
We refer interested readers to \cite{Asmussen2003} for a comprehensive overview of ladder processes.
For any non-decreasing function $F:\R\to\R$, we let $\Vert F \Vert := \lim_{x \rightarrow \infty} F(x)$. A distribution, $F$, is called proper if $\Vert F \Vert = 1$ and defective if $\Vert F \Vert < 1$.
There is a well-known result that $G_+$ is proper and $G_-$ is defective since $\mu_Y > 0$ \citep[cf. Theorem VIII 2.4.][]{Asmussen2003}.
This implies that the probability that $\tau_-$ is finite cannot be one, i.e. $ \lim_{x \to \infty} \P \{\tau_- < x\} <1$, whereas $\tau_+ < \infty$ almost surely. Additionally, $\E[\tau_+] < \infty$, whereas $\E[\tau_-]$ is infinite.
By Wald's identity \citep[cf. Appendix A][]{Tijm2003}, $\mu_+$ can be expressed as a function of $\E[\tau_+]$ and $\mu_Y$:
\[
    \mu_+ = \E [S_{\tau_+}] = \E \left[\textstyle\sum_{t=1}^{\tau_+} Y_t \right] = \E [\tau_+] \E [Y] = \E [\tau_+] \mu_Y.
\]
Let $m_n$ denote the partial minimum of the random walk within the first $n$ periods, i.e. $m_n := \min_{0 \leq t < n} S_t$. Then, the minimum of the entire random random walk, $m$, is defined as
$m := \inf_{0 \leq t < \infty} S_t$.
Define the descending ladder height renewal measure $U_- (x) := \sum_{t=0}^{\infty} G_-^{*t} (x)$, where $G_-^{*t}$ denotes the $t$-fold convolution of $G_-$, i.e., $G_-^{*t+1} (x) := G_-^{*t} * G_- (x)$, and $G_-^{*0} (x) = \delta_0 (x)$, with $\delta_0$ representing the probability measure degenerate at 0, i.e. $\delta_0(x)=1$ if $x\geq0$ and zero otherwise. We can express the distribution function of $m$ as \citep[cf. Theorem VIII, 2.2. ][]{Asmussen2003}:
\begin{equation}
    \label{eq:mprob}
    \P \{ m \leq x\} = \frac{U_- (x)}{\Vert U_- \Vert}.
\end{equation}
Next, $J_t$ is distributed as the waiting time of the $t$-th customer of a GI/G/1 queue with inter-arrival distribution $F_D$ and service time $r$. Thus, similar to Proposition, X.1.1. of \cite{Asmussen2003} $\Jinf \eqd -m$ ($\eqd$ denotes equality in distribution)
% \[
%     \P (\Jinf \leq x) = 1- \P (m \leq -x), \qquad x \geq 0,
% \]
which implies by Equation \eqref{eq:mprob} that:
\begin{equation}
    \label{eq:EJinf}
    \E[\Jinf] = - \E[m] = \frac{1}{\Vert U_- \Vert} \int_{-\infty}^0 U_-(x) dx.
\end{equation}
Similar to $U_-$, we define the ascending ladder height renewal measure, $U_+$, by $U_+ := \sum_{t=0}^{\infty} G_+^{*t}$.

\subsection{Solution to the Wiener-Hopf equation}
\label{section:WHequation}
We build on the methodology developed by \cite{Asmussen1998} to derive a solution to Equation \eqref{eq:RelVal}. \cite{Asmussen1998} shows that a solution to the Wiener-Hopf equation satisfies $v_r(x) = a_r*U_-*U_+ (x)$. Using this fact leads, after multiple intricate steps, to the characterization of $v_r(x)$ in Theorem \ref{theorem:solution}:
%A key distinction between our approach and that of \cite{Asmussen1998} lies in the class of admissible solutions: While \cite{Asmussen1998} restricts attention to non-negative solutions, we allow for all possible solutions, including non-positive ones.
%We use the following lemma to solve Equation \eqref{eq:RelVal} under this general class of admissible solutions.

%\begin{lemma} \citep[Corollary 3.1 and Proposition 3.3 of][]{Asmussen1998}
%    \label{lemma:WHsolution}
%    \[
%        v_r(x) = a_r*U_-*U_+ (x), \qquad v_r(0)=0, \forall x\geq0.
%    \]
%\end{lemma}

%\noindent
%Lemma \ref{lemma:WHsolution} provides a powerful approach for solving the Wiener-Hopf equation \eqref{eq:RelVal}.
%Applying Lemma \ref{lemma:WHsolution} to derive $v_r (x)$ is intricate and involves multiple steps.  
%We present the solution to Equation \eqref{eq:RelVal} in the next theorem.  

\begin{theorem}
    \label{theorem:solution}
    The relative value function characterized by Equation \eqref{eq:RelVal} is given by
    % \[
    %     v_r(x) = h \E[\tau_+] \int_0^x U_+ (y) dy - (h+p) x, \qquad x \geq 0.
    % \]
    \[
        v_r(x) = \frac{h \mu_+}{\mu_D-r} \int_0^x U_+ (y) dy - (h+p) x, \qquad x \geq 0.
    \]
\end{theorem}

\noindent
Observe that the ascending ladder process is in fact a renewal process. As such, it possesses all the general properties of the renewal processes including the following lemma.
Let $\kappa\in\R_+$ be expressed by
\[
    \kappa:=
    \begin{dcases}
        \frac{\sigma_+^2+\mu_+^2}{2\mu_+^2} & \text{if } D \text{ is non-lattice},\\
        \frac{\sigma_+^2+\mu_+^2+\mu_+}{2\mu_+^2} &\text{if } D \text{ is lattice.}
    \end{dcases}
\]

\begin{lemma}
    \label{lemma:UplusFunc}
    The ascending ladder height renewal measure $U_+$ can be expressed as 
    \[
        U_+ (x) = \frac{1}{\mu_+}x+\kappa + g_r(x).
    \]
    where $g_r:\R_+ \to \R$ satisfies $\vert g_r(x) \vert \leq \kappa$ for all $x\geq0$, and $g_r(x)=o(1)$.
\end{lemma}

\cite{vanJaarsveld2024} show that $v_r$ is a quadratic function in the case of exponential demand.
Next, we demonstrate that for a general demand distribution, $v_r$ can be expressed as the sum of a quadratic function and an $o(x)$ term. This result holds under the sole mild assumption that the demand distribution has a finite second moment.

\begin{theorem}
    \label{theorem:valuefunctionExpression}
    If $r \in [0,\mu_D)$ then for all $x\geq 0$,
    \[
        v_r(x) = b(r) \left(\left(x-\xiT(r)\right)^2 -\xiT^2(r) + 2\mu_+  \int_0^x g_r(y) dy \right), 
    \]
    with,
    \[
        b(r) := \frac{h}{2(\mu_D-r)}, \quad  \xiT(r):= (\mu_D-r) \left(\frac{p}{h} +1 \right)-\kappa 
    \]
    and $g_r$ as specified in Lemma \ref{lemma:UplusFunc}.
\end{theorem}
\proof{Proof.} This follows after some computation from Theorem \ref{theorem:solution} and Lemma \ref{lemma:UplusFunc}.
    \Halmos
\endproof

\noindent
Indeed, $g_r$ vanishes faster than $o(1)$ for most practical demand processes. For instance, it decays exponentially fast, i.e., $g_r(x) = o(e^{-\alpha x})$ with $\alpha>0$, if $D$ is non-lattice and sub-exponential, i.e., $\int_0^{\infty} e^{\delta x} G_+(dx) < \infty$, for some $\delta >0$ \citep[cf. VII Section 2.][]{Asmussen2003}. In this case, $v_r(x)$ is asymptotically quadratic as $x$ grows large.

% Although $\lim_{x \to \infty}\int_0^x g_r(y)dy$ may diverge, for any fixed $\vert \Delta x \vert < \infty$, the integral $\int_x^{x+\Delta x} g_r(y)dy$ vanishes as $x\to \infty$.
% It follows that for every $\Delta x \in \R$ and $\epsilon>0$ there exists a $x_{\epsilon} \geq 0$ such that:
% \begin{equation}
%     \label{eq:DifferenceValueFunc}
%     \vert v_r(x + \Delta x) - v_r(x) - b(r) \left( \left(x+ \Delta x - \left(\xi(r)-r\right)\right)^2-\left(x - \left(\xi(r)-r\right)\right)^2 \right) \vert \leq \epsilon, \qquad \forall x \geq x_{\epsilon}.
% \end{equation}
Next, suppose that $Z$ is a non-negative random variable. We introduce a sufficient condition for the existence and finiteness of $\E[v_r(Z)]$. 
\begin{lemma}
    \label{lemma:Ev_finiteness}
    Let $Z$ have a finite second moment and $r \in [0, \mu_D)$, then $\vert \E[v_r(Z)] \vert <\infty$.
\end{lemma}

\subsection{Inventory dynamics}
\label{section:InventoryDynamicsAsymp}
We next establish useful properties of constant order policies and PIL policies.
Recall that $r_p \in [0,\mu_D)$ represents the best constant order quantity under a lost sales unit penalty cost of $p \in \R_+$, given a fixed holding cost $h$, i.e., $r_p \in \arg\min_{r \in [0,\mu_D)} \C(C_r \mid p,L)$.  
As $p$ increases, it is intuitive that $r_p$ converges to $\mu_D$ to minimize the expected lost sales cost.  
In this case, the steady-state inventory level $J_{\infty} (C_{r_p})$ grows large as $\mu_Y \to 0$.
Next, we provide a more detailed elaboration on this intuition.
Consider the sequences of the steady state inventory levels $\{\Jinf(C_r)\}_{r \in [0,\mu_D)}$ .
\begin{lemma}
    \label{lemma:COPheavyTraffic}
    % \citep[Application of Theorem X 7.1][]{Asmussen2003}
    \begin{enumerate}[(a)]
        \item \label{lemma:Em_a} $\E[\Jinf (C_r)]$ is non-decreasing and convex in $r$,
        \item \label{lemma:Em_c} $2\mu_Y \E[\Jinf(C_r)]/\sigma^2_D \to 1$ as  $r \to \mu_D$.
        \item \label{lemma:rStar_monotonicity} $r_p$ is non-decreasing in $p\geq 0$,
        \item \label{lemma:rStar_limit}$r_p \to \mu_D$, as $p\to \infty$,
        \item \label{lemma:rStar_rate}$
        \sqrt{\frac{2p}{\sigma_D^2 h}}(\mu_D - r_p) \to 1$, as $p\to \infty$.
    \end{enumerate}
\end{lemma}
Note that part \ref{lemma:rStar_rate} of Lemma \ref{lemma:COPheavyTraffic} implies that $\lim_{p\to\infty}(r_p-\mu_D)/\sqrt{p}\in(0,\infty)$, that is, $r_p$ approaches $\mu_D$ at the same rate as $1/\sqrt{p}$ approaches 0.
\cite{Bu2020} have also established the same result in a more general setting with stochastic supply yield. Our setting allows for a shorter proof that we provide to make our argument self-contained.

We now shift our attention to the dynamics of the inventory level under PIL policies.
Let $\{\{q_t(P_{\xi})\}_{t\in \N_0}\}_{\xi\geq0}$ be a sequence of random variables representing orders under PIL policies $\{P_{\xi}\}_{\xi \geq 0}$, where $\{q_t(P_{\xi})\}_{t=0,\dots,L-1}$ are fixed for all $\xi\geq 0$ and known almost surely.
Let $\{\{J_t(P_{\xi})\}_{t\in \N_0}\}_{\xi\geq0}$ be the corresponding sequence of inventory level random variables.

% \begin{lemma}
%     \label{lemma:Jicv}
%     If $0 \leq \xi_1 \leq \xi_2$ then $J_t(P_{\xi_1})$ is less than $J_t(P_{\xi_2})$ in the increasing concave order, i.e. $J_t(P_{\xi_1}) \leq_{\text{icv}} J_t(P_{\xi_2})$ for all $t\in \N_0$.
% \end{lemma} 
% Lemma \ref{lemma:Jicv} shows that $J_t(P_{\xi})$ is both non-decreasing and concave in $\xi$ for any fixed demand realization.
\begin{lemma}
    \label{lemma:qbounded}
    For all $t\geq L+1$ and $\xi \geq 0$,
    the order size $q_t(P_{\xi})$ satisfy: $\E[q_t(P_{\xi}) ]\leq \min\{\xi, \mu_D\}$.
\end{lemma}
% Lemma \ref{lemma:qbounded} provides a uniform upper bound on $\E[q_t(P_{\xi})]$ and a uniform lower bound on $\E[J_t(P_{\xi})]$ for all $t\geq L+1$.
% % , which will be instrumental in the analysis that follows.
% In particular, it implies that $\E[J_t(P_{\xi})]$ grows linearly with $\xi$ as $\xi \to \infty$.
The result from Lemma \ref{lemma:qbounded} leads to Lemma \ref{lemma:deltaEVJqAsymtote} which will later be used to show that the impact of $g_r$ on the optimality gap of the PIL policy remains bounded for any $p$ and large $L$.
This is the sense in which the term of $v_r(x)$ that grows sublinearly becomes unimportant as $p$ grows.
\begin{lemma}
    \label{lemma:deltaEVJqAsymtote}
    There exists $0 \leq M<\infty$ such that for all $ t\geq L+1$, $r \in [0,\mu_D)$, and $\xi \geq 0$,
    \[
        b(r)\mu_+\E \left[\int_{J_t\left(P_{\xi}\right)}^{J_t\left(P_{\xi}\right)+q_t(P_{\xi})-r} g_r(y) dy \right] < M.
    \]
    % \[
    %     \E \left[\int_{J_t\left(P_{\xi(r_p)}\right)}^{J_t\left(P_{\xi(r_p)}\right)+\Delta q} g_r(y) dy \vert \xb_{t-L} \right] < \epsilon.
    % \]
\end{lemma}
% The following lemma strengthens this result by showing that the probability of $J_t(P_{\xi})$ remaining small becomes negligible as $\xi$ grows large.

% \begin{lemma}
%     \label{lemma:Junbounded}
%     For any $\epsilon \geq 0$ and $x \geq 0$, there exists $\xi_{x,\epsilon} \geq 0$ such that for all $t \geq L$, and $\xi \geq \xi_{x,\epsilon}$, $\P\{J_t(P_{\xi}) \leq x\} \leq \epsilon$.
% \end{lemma}

Next, we define the projected inventory level $\xi(r)$, $r \in [0,\mu_D)$ by
\begin{equation}
    \label{eq:xir}
    \xi(r):= \xiT(r) + r = \frac{\mu_Yp}{h}+ \mu_D - \kappa.
\end{equation}
Notice that, by Lemma \ref{lemma:COPheavyTraffic}\ref{lemma:rStar_rate}, in combination with Lemma \ref{lemma:momentsPlus} and Equation \eqref{eq:xir}, we have
\begin{equation}
    \label{eq:xiAsymptGrowth}
    0 < \xi(r_p)/\sqrt{p} < \infty, \qquad \text{as } p \to \infty.
\end{equation}
That is, $\xi(r_p)$ goes to infinity in the order of $\sqrt{p}$, as $p \to \infty$.
The next section provides a more detailed analysis of the cost-rate difference between $C_{r_p}$ and $P_{\xi(r_p)}$.

\subsection{Cost-rate difference between PIL and constant order policy}
\label{section:additionalCostRate}
Next, we derive an upper bound on the cost rate of a family of PIL policies by comparing it to that of corresponding constant order policies. 
%We then analyze this bound in the asymptotic regime as $p \to \infty$, leveraging the results from the previous section. This analysis leads to the proof of Theorem \ref{theorem:epsilonOptimality}.
One classical way of comparing the performance of two policies is by using the improvement theorem \citep[cf. Theorem 6.2.1.][]{Tijm2003}.
In general, applying the improvement theorem to our problem requires the consideration of $L+1$-dimensional state space.
However, the state space collapses to one dimensional for a system under a constant order policy since all order quantities as identical.
Lemma \ref{lemma:ImprovementLemma} adapts the improvement theorem for a constant order policy.
\begin{lemma}
    \label{lemma:ImprovementLemma}
    \citep[Similar to Lemma 4 of][]{vanJaarsveld2024}
    Let $t_1 \leq t_2$, $t_1,t_2 \in \N$ and suppose $q_t=r$ for all $t \in \{t_1, \dots, t_2\}$. Then,
    \[
        \E[c_{[t_1,t_2]} (C_r) \vert J_{t_1}] = v_r(J_{t_1}) - \E[v_r(J_{t_2+1}) \vert J_{t_1}] + (t_2+1-t_1) \C(C_r).
    \]
\end{lemma}

We are now in the position to prove the main results. 
% Let $\xi_p \in \arg \min_{\xi\geq 0} \C(P_{\xi})$ for $p\geq 0$.

\begin{theorem}
    \label{theorem:AsynmptoticOptimality}
    There exists $0 \leq M <\infty$ such that $\C(P_{\xi_p}) \leq \C(P_{\xi(r_p)}) \leq \C(C_{r_p}) + M$ for every $p \geq 0$. 
\end{theorem}

\proof{Proof of Theorem \ref{theorem:AsynmptoticOptimality}.}
    In this proof we bound $\E[c[L,T]\left(P_{\xi(r_p)}\right)- c[L,T]\left(C_{r_p}\right)]$ for $p \geq 0$. Similar to \cite{vanJaarsveld2024}, the proof relies on a policy $\PC^{\tT}$, $\tT\in \N_0$, which uses the PIL policy $P_{\xiT(r_p)}$ to order for $t=1, \dots, \tT+L$, and then orders $r_p$ when $t\geq \tT+L+1$, that is,
    \[
        \PC^{\tT}_t(\xb) = \begin{cases}
            P_{\xi(r_p)} (\xb), & \qquad t=1, \dots, \tT+L,\\
            r_p, & \qquad t=\tT+L+1, \dots
        \end{cases}
    \]
    Then,
    \begin{alignat}{2}
    \label{eq:Thm2proof1}
        & \E[c_{[L, T]} (\PC^{\tT}) - c_{[L, T]} (\PC^{\tT-1})]= 
        \E[c_{[\tT+L, T]} (\PC^{\tT}) - c_{[\tT+L, T]} (\PC^{\tT-1})]= \nonumber\\
        & \E\left[\E \left[c_{[\tT+L, T]} (C(r_p)) \vert  J_{\tT+L} = J_{\tT+L}(\PC^{\tT})+ q_{\tT+L}(\PC^{\tT})-r_p\right] - \E\left[c_{[\tT+L, T]} (C(r_p)) \vert  J_{\tT+L} = J_{\tT+L}(\PC^{\tT}) \right]\right].
    \end{alignat}
    The first equality in \eqref{eq:Thm2proof1} holds because $c_t$ remains the same under $\PC^{\tT-1}$ and $\PC^{\tT}$ for $t\leq \tT+L-1$.
    To justify the second equality, first observe that $J_{\tT+L-1}$ remains unchanged under $\PC^{\tT-1}$ and $\PC^{\tT}$ due to the dynamics of the inventory levels.
    Second, observe that under both policies the system receives $r_p$ in periods $t>\tT+L$.
    Third, observe that a system initiated at $ J_{\tT+L}(\PC^{\tT})$ and receiving the order quantity $q_{\tT+L}(\PC^{\tT})$ is equivalent to one starting at $J_{\tT+L}(\PC^{\tT})+ q_{\tT+L}(\PC^{\tT})-r_p$ and receiving an order quantity $r_p$. Thus, the second equality compares the total cost of two systems under $C_{r_p}$ with different initial inventory levels.
    Next, we use Lemma \ref{lemma:ImprovementLemma} to expand Equation \eqref{eq:Thm2proof1} as follows
    \begin{multline*}
        \E\left[c_{[L, T]} (\PC^{\tT}) - c_{[L, T]} (\PC^{\tT-1}) \right] =
        \E \Big[v_{r_p}\left(J_{\tT+L}(\PC^{\tT}) + q_{\tT+L}(\PC^{\tT}) -r_p \right) - v_{r_p}\left(J_{\tT+L}(\PC^{\tT+1})\right) \\
        - v_{r_p}\left(J_{T+1}(\PC^{\tT})\right)+v_{r_p}\left(J_{T+1}(\PC^{\tT})\right) \Big].
    \end{multline*}

    \noindent
    Let $\PC^{-1}:=C_{r_p}$. We use a telescopic sum, 
    \begin{alignat}{2}
        \label{eq:telescopicSum}
        &\E \left[c_{[L,T]}\left(P_{\xi(r_p)} \right) - c_{[L,T]} \left(C_{r_p} \right) \right] = \E \left[c_{[L, T]} \left(\PC^{T} \right) - c_{[L, T]} \left(\PC^{-1} \right) \right]= \sum_{\tT=0}^{T-L} \E[c_{[L,T]}(\PC^{\tT}) - c_{[L,T]}(\PC^{\tT-1})]\nonumber\\
        & = \sum_{\tT=0}^{T-L} \E[v_{r_p}(J_{\tT+L}(\PC^{\tT}) + q_{\tT+L}(\PC^{\tT}) -r_p) - v_{r_p}(J_{\tT+L}(\PC^{\tT}))]  - \E[v_{r_p}(J_{T+1}(\PC^{\tT}))]+\E[v_{r_p}(J_{T+1}(\PC^{\tT-1}))].
    \end{alignat}
    By Theorem \ref{theorem:valuefunctionExpression} we notice that for any $r\in [0,\mu_D)$, $\xb_{t}\in \R_+^{L+1}$, $t\geq 0$, $\Tilde{q} \in \R$
    \begin{alignat}{2}
        \label{eq:EVr}
        & \E \left[ v_r(J_{t+L}+\Tilde{q}) \vert \xb_{t} \right] =  b(r) \E \left[\left(J_{t+L}+\Tilde{q}-\xiT(r)\right)^2 -\xiT^2(r) + 2\mu_+ \int_0^{J_{t+L}+\Tilde{q}} g_r(y) dy \vert \xb_{t} \right] = \nonumber\\
        & b(r) \left(\Var[J_{t+L} \vert \xb_{t}]+ \left( \E[J_{t+L} \vert \xb_{t}] +\Tilde{q} - \xiT(r) \right)^2 -\xiT^2(r) + 2\mu_+ \E \left[ \int_0^{J_{t+L}+\Tilde{q}} g_r(y) dy \vert \xb_{t} \right] \right).
    \end{alignat}
    Using Equation \eqref{eq:telescopicSum} combined with \eqref{eq:EVr} and some algebra we have
    \begin{multline*}
        % \label{eq:CostDifffinitTime}
        \E[c_{[L,T]}(P_{\xi(r_p)}) - c_{[L,T]}(C_{r_p})] =   \E[v_{r_p}(J_{T+1}(P_{\xi (r_p)}))] - \E[v_{r_p}(J_{T+1}(C_{r_p}))] +\\
        \sum_{\tT=0}^{T-L} b(r_p)\E \left[ -(P_{\xi(r_p)}(\xb_{\tT})-r_p)^2 + 2\mu_+ \E \left[\int_{J_{\tT+L}(\PC^{\tT})}^{J_{\tT+L}(\PC^{\tT})+P_{\xi(r_p)}(\xb_{\tT})-r_p} g_{r_p}(y) dy \Big\vert \xb_{\tT} \right] \right].
    \end{multline*}
    % Now by Lemma \ref{lemma:deltaEVJqAsymtote}, for any $\epsilon >0$ there exists a $p_{\epsilon}$ such that for all $p \geq p_{\epsilon}$ and $\tT \in \N_0$:
    % \[
    %     \E \left[\int_{J_{\tT+L}(\PC^{\tT})}^{J_{\tT+L}(\PC^{\tT})+P_{\xi(r_p)}(\xb_{\tT})-r_p} g_{r_p}(y) dy \Big\vert \xb_{\tT} \right] = \E \left[\int_{J_{\tT+L}(P_{\xi(r_p)})}^{J_{\tT+L}(P_{\xi(r_p)})+P_{\xi(r_p)}(\xb_{\tT})-r_p} g_{r_p}(y) dy \Big\vert \xb_{\tT} \right] < \epsilon.
    % \]
    Then, it follows from Lemma \ref{lemma:deltaEVJqAsymtote} that there exists $0 \leq M < \infty$ such that for all $p \geq 0$ 
    \begin{alignat}{2}
        \label{eq:PILFiniteBound}
        &\E[c_{[L,T]}(P_{\xi(r_p)}) - c_{[L,T]}(C_{r_p})] < \nonumber\\
        &\E \left[v_{r_p}(J_{T+1}(P_{\xi(r_p)})) \right] - \E \left[v_{r_p}(J_{T+1}(C_{r_p})) \right] + (T-L+1) M
         - b(r_p)\sum_{\tT=0}^{T-L}\E \left[(P_{\xi(r_p)}(\xb_{\tT})-r_p)^2 \right] \nonumber\\
         &\leq \E \left[v_{r_p}(J_{T+1}(P_{\xi(r_p)}))] - \E[v_{r_p}(J_{T+1}(C_{r_p})) \right] + (T-L+1) M.
    \end{alignat}
    Notice that the last inequality of \eqref{eq:PILFiniteBound} holds since $b(r_p)$ and $(P_{\xi(r_p)}(\xb_{\tT})-r_p)^2$ are non-negative.
    We use \eqref{eq:PILFiniteBound} to find an upper bound on the cost-rate of the PIL policy, i.e., $\C(P_{\xi(r_p)})$ with respect to the cost-rate of constant order policy $\C(C_{r_p})$ when $p\geq 0$:
    \begin{alignat}{2}
        \label{eq:PILInfiniteBound}
        C(P_{\xi(r_p)}) =& \limsup_{T \to \infty} \frac{1}{T-L+1} \E[c_{[L,T]}(P_{\xi(r_p)})] \nonumber\\
        <& \limsup_{T \to \infty} \frac{1}{T-L+1} \left(\E[c_{[L,T]}(C_{r_p}) + v_{r_p}(J_{T+1}(P_{\xi(r_p)})) - v_{r_p}(J_{T+1}(C_{r_p}))] + (T-L+1) M \right) \nonumber\\
        =& C(C_{r_p}) + M + \limsup_{T \to \infty} \frac{1}{T-L+1} \left(\E[v_{r_p}(J_{T+1}(P_{\xi(r_p)}))] - \E\left[v_{r_p}(J_{T+1}(C_{r_p}))\right]\right) \nonumber\\
        =& C(C_{r_p}) + M
    \end{alignat}
    The last equality holds since both $\E[v_{r_p}(J_{T+1}(P_{\xi(r_p)}))$ and $\E[v_{r_p}(J_{T+1}(C_{r_p}))]$ remain finite as $T \to \infty$.
    First, observe, as \cite{vanJaarsveld2024} do, that $0 \leq J_{T+1}(P_{\xi(r_p)}) \leq \xi(r_p) + L \mu_D$ for all $T \geq L$ which ensures that $J_{T+1}(P_{\xi(r_p)})$ has finite first and second moments. %\cite{vanJaarsveld2024} provide an interesting argument on the upper-bound of $J_{T+1}(P_{\xi(r_p)})$. 
    Then it follows from Lemma \ref{lemma:Ev_finiteness} that $\vert \E[v_{r_p}(J_{T+1}(P_{\xi(r_p)})) \vert < \infty$.
    Second, $J_{T+1}(C_{r_p})$ converges to the steady state distribution of the inventory level under the constant order policy $C_{r_p}$, i.e., $J_{\infty}$, as $T \to \infty$. We note that $J_{\infty}$ has a finite first moment because $D$ has a finite second moment. Additionally, $J_{\infty}$ has a finite second moment since $0 \leq \left((r_p- D)^+\right)^3 \leq r_p^3$, implying that $\E \left [\left((r_p- D)^+ \right)^3 \right] < \infty$ \citep[cf. Theorem X. 2.1.][]{Asmussen2003}. Hence, $\vert \E[v_{r_p}(J_{T+1}(C_{r_p}))] \vert < \infty$ due to Lemma \ref{lemma:Ev_finiteness}.
    The optimality of $\xi_p$, i.e., $C(P_{\xi_p}) \leq C(P_{\xi(r_p)})$ together with Inequality \eqref{eq:PILInfiniteBound} complete the proof.
    \Halmos
\endproof

\proof{Proof of Theorem \ref{theorem:epsilonOptimality}.}
    Combining Theorem \ref{theorem:AsynmptoticOptimality} with asymptotic optimality of the constant order policy as $L$ approaches infinity \citep{Goldberg2016,Xin2016} provides the result:
    There exists $0 \leq M < \infty$ such that for all $p \geq0$:
    \[
        \lim_{L \to \infty} \left(\C(P_{\xi_p}) - \C^*(p,L) \right) \leq \lim_{L \to \infty} \left(\C(P_{\xi(r_p)}) - \C^* (p,L) \right) \leq M.
    \]
    Observe that $\lim_{p\to \infty} \lim_{L \to \infty} \C^*(p,L) = \infty$.
    Thus,
    \[
        \lim_{p\to \infty} \lim_{L \to \infty} \frac{\C(P_{\xi_p})}{\C^*(p,L)} = 1. \Halmos
    \]
\endproof
% \subsection{Asymptotic Optimality as $L \to \infty$}
% We combine Theorem \ref{theorem:AsynmptoticOptimality} with asymptotic optimality of the constant order policy as $L$ approaches infinity \citep{Goldberg2016,Xin2016} to introduce the following asymptotic optimality result for the PIL policy. Let $\pi^*$ denote the optimal policy.

% \begin{theorem}
%     \label{theorem:epsilonOptimality}
%     For any $\epsilon >0$, there exists $p_{\epsilon} \geq 0$ such that the PIL policy is asymptotically $\epsilon-$optimal for every $p \geq p_{\epsilon}$ when demand has a finite second moment:
%     \[
%         \lim_{L \to \infty} C(P_{\xi(r_p)}) - C(\pi^*) < \epsilon, \qquad p \geq p_{\epsilon}.
%     \]
% \end{theorem}
% It is worth mentioning that the cost-rate of the optimal policy tends to infinity as $p$ grows large.
% Theorem \ref{theorem:epsilonOptimality} implies that under such situation the cost-rate of the PIL policy remains in the $\epsilon-$neighborhood of the optimal cost-rate if $L \to \infty$.

\section{Concluding remarks}
\label{Section:Conclusion}
In this paper, we proved that the PIL policy is asymptotically optimal for sufficiently large lost sales unit costs as the leadtime approaches infinity, under mild assumptions on the i.i.d. demand process. This result, combined with \cite{vanJaarsveld2024}, demonstrates that the PIL policy is asymptotically optimal when the lost sales penalty cost is large, both in the case of a small leadtime and when the leadtime grows at a rate faster than the unit cost of lost sales. This makes the PIL policy the only single-parameter policy that guarantees optimality in both regimes under a general i.i.d. demand.
% We also demonstrated that the optimal cost-rate approaches infinity proportional to the square root of the lost sales unit penalty cost when both leadtime and lost sales unit penalty cost approach infinity with the leadtime growing at a faster rate.
It remains an open question whether the PIL policy is asymptotically optimal when both the leadtime and the lost sales unit penalty cost grow at the same rate. To the best of our knowledge, no simple policies are known to achieve optimality in this asymptotic regime.

\begin{APPENDIX}{}
\section{Proof of Lemma \ref{lemma:momentsPlus}}
\label{section: ProofmomentsPlus}
Part (a) follows from Theorem XVIII.5.1. \cite{Feller1991V2}.
The rest is the proof of Part (b).
Observe that $\E[Y^2]<\infty$ only if for some $\alpha>2$,
\begin{equation*}
    \label{eq:regVarCond}
    1-F_Y(x) = O(x^{-\alpha}) \quad\text{as } x \to \infty.
\end{equation*}
This condition is equivalent to
\begin{equation}
    \label{eq:partialMoments}
     \lim_{x \to \infty} \E\left[Y^2 \vert Y \geq x\right]=\lim_{x \to \infty} \frac{\int_x^{\infty} y^2 F_Y(dy)}{1-F_Y(x)} = \lim_{x \to \infty} \frac{x^2(1-F_Y(x))+2 \int_x^{\infty}y(1-F_Y(y))dy}{1-F_Y(x)} < \infty.
\end{equation}
Since for any $x \in (0,\infty)$, $\E\left[Y^2 \vert Y \geq x\right] < \infty$, and it is finite at the limit $x \to \infty$ by \eqref{eq:partialMoments}, we conclude that,
\begin{equation}
    \label{eq:partialMomentsSup}
    \sup_{x\geq0} \E\left[Y^2 \vert Y \geq x\right]< \infty.
\end{equation}
Next, we notice that $S_{\tau_+} = S_{\tau_+-1}+Y_{\tau_+} \eqd S_{\tau_+-1}+Y \vert Y\geq -S_{\tau_+-1}$, and $S_{\tau_+-1}<0$ by the definition of $\tau_+$. This implies in particular that $S_{\tau_+}\leq Y_{\tau_+}$ almost surely and 
\[
    \E\left[S_{\tau_+}^2\right] \leq \E[Y_{\tau_+}^2] = \E\left[\E\left[Y^2 \vert Y\geq -S_{\tau_+-1}\right]\right] \leq \sup_{x>0} \E\left[Y^2 \vert Y \geq x\right] < \infty.
\]
Finally $\E[D^2] < \infty$ is equivalent to $\E[Y^2] <\infty$ which completes the proof.
\Halmos

\section{Proof of Theorem \ref{theorem:solution}}
\label{section:ProofThmSolution}
We use the methodology of solving Wiener-Hopf equations introduced by \cite{Asmussen1998}.
A key distinction between our approach and that of \cite{Asmussen1998} lies in the class of admissible solutions: While \cite{Asmussen1998} restricts attention to non-negative solutions, we allow for all possible solutions, including non-positive ones.
We use the following lemma to solve Equation \eqref{eq:RelVal} under this general class of admissible solutions.
\begin{lemma} \citep[Corollary 3.1 and Proposition 3.3 of][]{Asmussen1998}
   \label{lemma:WHsolution}
   \[
       v_r(x) = a_r*U_-*U_+ (x), \qquad v_r(0)=0, \forall x\geq0.
   \]
\end{lemma}

\noindent
Lemma \ref{lemma:WHsolution} provides a powerful approach for solving the Wiener-Hopf equation \eqref{eq:RelVal}.
Applying Lemma \ref{lemma:WHsolution} to derive $v_r (x)$ is intricate and involves multiple steps.
The proof of Theorem \ref{theorem:solution} is provided at the end of this section.
The first step in deriving $v_r(x)$, following Lemma \ref{lemma:WHsolution}, involves expressing $a_r (x)$ in terms of $F_Y (x)$.
This step is necessary due to the lack of a standard result in the literature that allows direct convolution of $a_r$ in Equation \eqref{eq:ar} with $U_-$.
However, as we will later show, existing results from random walk theory enable the convolution of $F_Y$ with both $U_-$ and $U_+$.

\begin{lemma}
    \label{lemma:ar}
    \[
        a_r(x) = (h+p) \int_{-r}^x  F_Y (y) dy - p x - h \E[\Jinf].
    \]
\end{lemma}

\proof{Proof of Lemma \ref{lemma:ar}.}
    By Equation \eqref{eq:ar}:
    \begin{align*}
        a_r(x) & = h \E_Y[(x-Y)^+] + p \E_Y[(Y-x)^+]- p \mu_Y + h \E[\Jinf] \\
        & = (h+p) \E_Y[(x-Y)^+] - p x + p \mu_Y - p \mu_Y + h \E[\Jinf].
    \end{align*}
    Now, we express $\E_Y[(x-Y)^+]$ in terms of $F_Y (x)$ as follows:
    \[
        \E_Y[(x-Y)^+] = \int_{-r}^x (x-y) F_Y(dy) =  x F_Y (y) \Big\vert_{-r}^x -\int_{-r}^x y F_Y(dy).
    \]
    By assumptions, $F_Y(-r) = 0$. Additionally, we use integration by parts to compute $\int_{-r}^x y F_Y(dy)$:
    \[
        \quad \int_{-r}^x y F_Y(dy) = y F_Y (y) \Big\vert_{-r}^x - \int_{-r}^x  F_Y (y) dy.
    \]
    Thus
    \[
        \E_Y[(x-Y)^+]=  \int_{-r}^x  F_Y (y) dy,
    \]
    and
    \[
        a_r(x) = (h+p) \int_{-r}^x  F_Y (y) dy - p x - h \E[\Jinf]. \Halmos
    \]
\endproof

\noindent
By Lemma \ref{lemma:ar}, $a_r$ is expressed as a linear combination of $\int_{-r}^x  F_Y (y) dy$, $x$, and the constant 1.
Importantly, the convolution operator possesses both distributive and homogeneous properties.
These properties enable the separate convolution of $\int_{-r}^x F_Y(y) dy$, $x$, and $1$ with $U_-$ and $U_+$, providing the basis for the proof of Theorem \ref{theorem:solution}.

\noindent\textbf{Derivation of} $\int_{-\infty}^x  F_Y (y) dy * U_- * U_+ (x)$\textbf{:}
The convolution operator satisfies the associativity property. Furthermore, the following well-known lemma indicates the relation between the integration and convolution operators.

% \begin{lemma}
%     \label{lemma:integral}
%     For two real functions $K$ and $N$, let $K(x) = \int_{-\infty}^x k(y) dy$, and $N(x) = \int_{-\infty}^x n(y) dy$. then,
%     \[
%         K*N(x) = N*K(x) = \int_{-\infty}^x k*N(y) dy =   \int_{-\infty}^x n*K(y) dy.
%     \]
% \end{lemma}

\noindent
Associativity and commutativity of convolution imply that:

\begin{equation}
    \label{eq:integralrelation}
    \int_{-\infty}^x  F_Y (y) dy * U_- * U_+ (x) =
    \int_{-\infty}^x F_Y * U_- (y) dy  * U_+ (x).
\end{equation}

\noindent
By Equation \eqref{eq:integralrelation}, the next steps involve first calculating $F_Y * U_-$, then convolving the result with $U_+$, and finally integrating the outcome. The following lemma is crucial to our computations.

\begin{lemma} \citep[Theorem VIII 3.1. and Corollary 3.2][]{Asmussen2003}
    \label{lemma:FUConvolution}
        \[
            U_-*F_Y = U_- + G_+ - \delta_0,
        \]
\end{lemma}

\noindent
It follows from commutativity of convolution and Lemma \ref{lemma:FUConvolution} that:
\begin{alignat}{2}
    \label{eq:FU1}
    & \int_{-\infty}^x U_-*F_Y(y)dy *U_+ (x) = \int_{-\infty}^x (U_- + G_+ - \delta_0)(y)dy *U_+ (x) = \nonumber\\
    & \int_{-\infty}^x (G_+ - \delta_0)*U_+(y)dy +  \int_{-\infty}^x U_- (y) dy * U_+ (x).
\end{alignat}
By definition of $U_+$:
\begin{equation}
    \label{eq:identityG+U}
    G_+ * U_+ = G_+ *  \sum_{t=0}^{\infty} G_+^{*t}= \sum_{t=1}^{\infty} G_+^{*t} = U_+ - \delta_0.
\end{equation}
Furthermore, it is a well-known result that the convolution of any function with $\delta_0$ returns the same function.
Consequently, the first term of Equation \eqref{eq:FU1} can be calculated as follows:
\begin{equation}
    \label{eq:FU1First}
    \int_{-\infty}^x (G_+ - \delta_0)*U_+(y)dy =
    \int_{-\infty}^x (U_+ - \delta_0 - U_+) dy =
    - \int_{0}^x \delta_0 dy = - x.
\end{equation}
Now we address the second term of Equation \eqref{eq:FU1},
\[
    \int_{-\infty}^x U_- (y) dy * U_+ (x) = \int_{-\infty}^0 U_- (y) dy * U_+ (x) +
    \int_0^x U_- * U_+ (y) dy.
\]
We notice that for all $x \geq 0$, $U_-(x) = \Vert U_- \Vert$. Additionally, $1 * U_+ =U_+$, since for all $x \leq 0$, $U_+(x) = 0$. 
Thus, by Equation \eqref{eq:EJinf}:
\begin{equation}
    \label{eq:FU1second1}
    \int_{-\infty}^x U_- (y) dy * U_+ (x) = \Vert U_- \Vert \E[\Jinf] U_+ + \Vert U_- \Vert \int_0^x U_+ (y) dy.
\end{equation}

\noindent
The following lemma allows us to relate Equation \eqref{eq:FU1second1} to $\E[\tau_+]$.
\begin{lemma} \citep[Theorem VIII 2.3. (c)][]{Asmussen2003}
    \label{lemma:etau+}
    \[
        \lVert U_- \rVert = \E [\tau_+] = (1- \Vert  G_- \Vert)^ {-1}.
    \]
\end{lemma}

\noindent
By Lemma \ref{lemma:etau+} and Equation \eqref{eq:FU1second1}, we can compute the second term of Equation \eqref{eq:FU1}:
\begin{equation}
    \label{eq:FU1secondfinal}
    \int_{-\infty}^x U_- (y) dy * U_+ (x) = \E [\tau_+] \E[\Jinf] U_+ + \E [\tau_+] \int_0^x U_+ (y) dy.
\end{equation}

\noindent
We combine Equations \eqref{eq:FU1}, \eqref{eq:FU1First}, and \eqref{eq:FU1secondfinal} to compute $\int_{-\infty}^x  F_Y (y) dy * U_- * U_+ (x)$:
\begin{equation}
    \label{eq:arU_U+first}
    \int_{-\infty}^x  F_Y (y) dy * U_- * U_+ (x) = \E [\tau_+] \int_0^x U_+ (y) dy + \E [\tau_+] \E[\Jinf] U_+ - x.
\end{equation}

\noindent\textbf{Derivation of} $x * U_- * U_+ (x)$\textbf{:}
It is straightforward to verify that $x * U_- (x) = \int_{-\infty}^x U_-(y)dy$, given the definition and commutativity of the convolution operator.
Thus, by Equation \eqref{eq:FU1secondfinal}:
\begin{equation}
    \label{eq:arU_U+second}
    x * U_- * U_+ (x) = \int_{-\infty}^x U_-(y)dy *U_+ (x) =  \E [\tau_+] \int_0^x U_+ (y) dy + \E [\tau_+] \E[\Jinf] U_+.
\end{equation}

\noindent\textbf{Derivation of} $1 * U_- * U_+ (x)$\textbf{:}
By definition of the convolution operator:
\[
    1 * U_- (x) = \int_{-\infty}^x U_-(dy) = U_-(x) = \Vert U_- \Vert = \E [\tau_+],
\]
which implies that:

\begin{equation}
    \label{eq:arU_U+third}
    1 * U_- * U_+ (x) = \E [\tau_+] U_+.
\end{equation}

\noindent
At this point we have all the tools available to prove Theorem \ref{theorem:solution}.

\proof{Proof of Theorem \ref{theorem:solution}.}
    By Lemma \ref{lemma:WHsolution} and Lemma \ref{lemma:ar}, for $x \geq 0$, $v_r (x)$ can be calculated by:
    \begin{alignat*}{2}
        & v_r (x) = a_r*U_-*U_+ (x) = \left ( (h+p) \int_{-r}^x  F_Y (y) dy - p x - h \E[\Jinf] \right)*U_-*U_+ (x),\\
        & = (h+p) \int_{-r}^x  F_Y (y) dy *U_-*U_+ (x) - p x *U_-* U_+ (x) - h \E[\Jinf] 1 * U_- * U_+ (x).
    \end{alignat*}

    \noindent
    By Equations \eqref{eq:arU_U+first}, \eqref{eq:arU_U+second}, and \eqref{eq:arU_U+third}:
    \begin{multline*}
        v_r (x) = (h+p) \left( \E [\tau_+] \int_0^x U_+ (y) dy + \E [\tau_+] \E[\Jinf] U_+ - x \right) + \\
        - p \left( \E [\tau_+] \int_0^x U_+ (y) dy + \E [\tau_+] \E[\Jinf] U_+ \right) - h \E[\Jinf] \E [\tau_+] U_+.
    \end{multline*}

    \noindent
    Simplifying the last expression, we can calculate $v_r (x)$ as follows:
    \[
        v_r (x) = h  \E [\tau_+] \int_0^x U_+ (y) dy - (h+p) x.
    \]
    By Wald's equality $\E[\tau_+] \mu_Y = \E[S_{\tau_+}] = \mu_+$, since $\tau_+$ is a stopping time for the $\{S_t\}_{t \in \N}$ process. Hence:
    \[
        v_r (x) = \frac{h \mu_+}{\mu_Y} \int_0^x U_+ (y) dy - (h+p) x.\Halmos
    \]
    
\endproof

\section{Proof of Lemma \ref{lemma:UplusFunc}}
\label{section:ProofUplusFunc}
For this proof we need two observations. First, for all $x \geq 0$:
\begin{equation}
    \label{eq:UplusFunc_bounds}
    \frac{1}{\mu_+}x \leq U_+ (x) \leq \frac{1}{\mu_+}x + \kappa.%\frac{\sigma_+^2+\mu_+^2}{\mu_+^2}.
\end{equation}
The left inequality of \eqref{eq:UplusFunc_bounds} deals with the fact that the expected time until the next renewal after $x$ (residual life) is non-negative \citep[cf. V. 6.][]{Asmussen2003}.
The right inequality of \eqref{eq:UplusFunc_bounds} is Lorden's Inequality \citep{Lorden1970}.
Next, as $x \to \infty$,
\begin{equation}
    \label{eq:UplusFunc_asymptote}U_+ (x) = \frac{1}{\mu_+}x+\kappa + o(1).
\end{equation}
Equation \eqref{eq:UplusFunc_asymptote} is due to the asymptotic expansion of the expected residual life function \citep[cf. Proposition V 6.1.][for non-lattice $D$]{Asmussen2003}.
\eqref{eq:UplusFunc_bounds} together with \eqref{eq:UplusFunc_asymptote} provide the result.
\Halmos

\section{Proof of Lemma \ref{lemma:Ev_finiteness}}
    \label{section:ProofEv_finiteness}
    By Theorem \ref{theorem:valuefunctionExpression},
    \[
        \left(Z-\xiT(r)\right)^2 - \kappa\leq \frac{1}{b(r)} \left(v_r(Z) + \xiT^2(r)\right) \leq \left(Z-\xiT(r)\right)^2 + \kappa.
    \]
    We take the expectation with respect to $Z$ on all sides,
    \[
        \E\left[\left(Z-\xiT(r)\right)^2 \right] - \kappa \leq \frac{1}{b(r)} \left(\E[v_r(Z)] + \xiT^2(r)\right) \leq \E\left[\left(Z-\xiT(r)\right)^2\right] + \kappa.
    \]
    Observe that by definition, $\Var \left[Z-\xiT(r) \right] = \E \left[\left(Z-\xiT(r)\right)^2 \right] - \left(\E[Z] - \xiT(r)\right)^2$.
    It follows that,
    \[
        \Var[Z]+\left(\E[Z]-\xiT(r)\right)^2 - \kappa\leq \frac{1}{b(r)} \left(\E[v_r(Z)] + \xiT^2(r)\right) \leq \Var[Z]+\left(\E[Z] - \xiT(r)\right)^2 + \kappa.
    \]
    Notice that $Z$ has finite first and second moments and $0 < \mu_+,\sigma_+ <\infty$ for $r \in [0,\mu_D)$, which implies the result.
    \Halmos

\section{Proof of Lemma \ref{lemma:COPheavyTraffic}}
\label{section:prooflemma:COPheavyTraffic}
Consider the sequences of random variables $\{\{Y_t(r) = D_t - r\}_{t \in \N}\}_{r \in [0,\mu_D)}$, sequences of random walks $\{\{S_t(r) = \sum_{i=1}^t Y_t(r)\}_{t \in \N}\}_{r \in [0,\mu_D)}$.
\begin{enumerate}[(a)]
    \item We recall that $\Jinf \eqd -m$. It is a known result \citep[cf. Proposition VIII 4.5][]{Asmussen2003} that,
    \[
        \E[\Jinf (C_r)] = \sum_{t=1}^{\infty} \frac{1}{t} \E[S_t^-] = \sum_{t=1}^{\infty} \frac{1}{t} \E[(-S_t)^+] = \sum_{t=1}^{\infty} \frac{1}{t} \E\left[\left(tr - \sum_{i=1}^t D_i \right)^+\right].
    \]
    Let $r_1, r_2 \in [0,\mu_D)$ and $r_1 \leq r_2$. First we prove monotonicity. Observe that 
    \[
        tr_1 - \sum_{i=1}^t D_i \leq tr_2 - \sum_{i=1}^t D_i,
    \]
    almost surely and so
    \[
        \E\left[(tr_1 - \sum_{i=1}^t D_i)^+ \right] \leq \E \left[(tr_2 - \sum_{i=1}^t D_i)^+ \right].
    \]
    Hence,
    \[
        \E[\Jinf (C_{r_1})] = \sum_{t=1}^{\infty} \frac{1}{t} \E\left[\left(tr_1 - \sum_{i=1}^t D_i \right)^+ \right] \leq  \sum_{t=1}^{\infty} \frac{1}{t} \E\left[\left(tr_2 - \sum_{i=1}^t D_i\right)^+ \right] = \E[\Jinf (C_{r_2})].
    \]
    Next, we prove convexity. For all $0 \leq \alpha \leq 1$,
    \begin{alignat*}{2}
        & \left(t(\alpha r_1 + (1-\alpha)r_2) - \sum_{i=1}^t D_i \right)^+ =
        \left(\alpha \left(t r_1 - \sum_{i=1}^t D_i\right) + (1-\alpha) \left(t r_2 - \sum_{i=1}^t D_i \right) \right)^+\\
        & \leq \alpha \left(t r_1 - \sum_{i=1}^t D_i \right)^+ + (1-\alpha) \left(t r_2 - \sum_{i=1}^t D_i\right)^+, \qquad \text{almost surely.}
    \end{alignat*}
    % By Theorem (1.A.1) \cite{shaked2007},
    % \[
    %     \left(t(\alpha r_1 + (1-\alpha)r_2) - \sum_{i=1}^t D_i \right)^+ \leq_{\text{st}} \alpha (t r_1 - \sum_{i=1}^t D_i)^+ + (1-\alpha) (t r_2 - \sum_{i=1}^t D_i)^+.
    % \]
    Hence,
    \[
        \E\left[\left(t(\alpha r_1 + (1-\alpha)r_2) - \sum_{i=1}^t D_i \right)^+ \right] \leq \alpha \E \left[\left(t r_1 - \sum_{i=1}^t D_i \right)^+ \right] + (1-\alpha) \E\left[\left(t r_2 - \sum_{i=1}^t D_i \right)^+ \right],
    \]
    which gives,
    \begin{alignat*}{2}
        & \E[\Jinf(C_{\alpha r_1 + (1-\alpha)r_2})] = \sum_{t=1}^{\infty}\E \left[\left(t(\alpha r_1 + (1-\alpha)r_2) - \sum_{i=1}^t D_i \right)^+\right]\\
        \leq & \alpha \sum_{t=1}^{\infty} \E\left[\left(t r_1 - \sum_{i=1}^t D_i \right)^+\right] + (1-\alpha) \sum_{t=1}^{\infty} \E \left[\left(t r_2 - \sum_{i=1}^t D_i \right)^+\right]
        = \alpha \E[\Jinf(r_1)] + (1-\alpha) \E[\Jinf(r_2)].
    \end{alignat*}
    \item Part \ref{lemma:Em_c} presents the expected waiting time of a GI/G/1 queue in a heavy traffic condition. Interested readers may refer to \cite{Kingman1961}.

\item Recall that $\C(C_r) = h\E[\Jinf(r)] + p (\mu_D - r)$.
By part \ref{lemma:Em_a}, $\C(C_r)$ is convex in $r$.
Let $\partial \C(C_r)$ denote the sub-differential of the cost-rate function at $r$, that is:
\[
    \partial \C(C_r) := \{x \in \R: \C(C_{\Bar{r}}) - \C(C_r) \geq x (\Bar{r} -r), \forall \Bar{r} \geq 0\}.
\]
By the optimality condition $0 \in \partial \C(C_{r_p})$ which is equivalent to $\frac{p}{h} \in \partial \E[\Jinf(C_{r_p})]$, $p \geq 0$.
It is straightforward to verify that $\partial \E[\Jinf(C_r)]$, $r \in[0,\mu_D)$ is an interval $[a_r,b_r]$ where $a_r,b_r$ are some non-negative real numbers due to part \ref{lemma:Em_a}. Additionally, for any $0 \leq r_1 \leq r_2 < \infty$, $b_{r_1} \leq a_{r_2}$ due to the convexity of $\E[\Jinf(C_r)]$.
This implies that for $p_1 \leq p_2$, $r_{p_1} \leq r_{p_2}$, since $p_i/h \in \partial\E \left[\Jinf \left(C_{r_{p_i}}\right) \right]$, for $i\in\{1,2\}$, and either $b_{r_{p_1}} \leq a_{r_{p_2}}$ or $b_{r_{p_2}} \leq a_{r_{p_1}}$. 

\item Next we prove that $r_p$ approaches $\mu_D$ as $p \to \infty$. This statement is equivalent to showing that there exists no $0 \leq \Tilde{r} < \mu_D$ such that for some $\Tilde{p} \geq0$, $r_p \leq \Tilde{r}$ for all $p \geq \Tilde{p}$, considering part \ref{lemma:rStar_monotonicity}.
Assume the contrary that there exist such $\Tilde{r}$ and $\Tilde{p}$. 
Consider the sequence $\{p_r = \max(\frac{h\sigma_D^2}{2(\mu_D-r)^2},\Tilde{p})\}_{r \in (\Tilde{r},\mu_D)}$.
By assumption and convexity of $\C(C_{r_p})$, for any $r \in (\Tilde{r}, \mu_D)$:
\[
    h \E[\Jinf(C_{\Tilde{r}})] + p_r (\mu_D - \Tilde{r}) \leq h \E[\Jinf(C_r)] + p_r (\mu_D - r).
\]
It follows that:
\[
    \frac{p_r(r-\Tilde{r})}{h} \leq \E[\Jinf(C_r)] - \E[\Jinf(C_{\Tilde{r}})] \leq \E[\Jinf(C_r)].
\]
Therefore by definition of $p_r$
\begin{equation}
    \label{eq:OptCondP_r}
    \frac{\sigma_D^2(r-\Tilde{r})}{2(\mu_D-r)^2} \leq \E[\Jinf(C_r)].
\end{equation}
Now by part \ref{lemma:Em_c}, for every $\epsilon >0$ there exists $\Tilde{r}_{\epsilon} \in [0,\mu_D)$ such that for all $r \geq \Tilde{r}_{\epsilon}$
\[
    \frac{\E[\Jinf(C_r)]}{\sigma_D^2/(\mu_D-r)} < 1+\epsilon.
\]
Let $\max\{\Tilde{r} , \Tilde{r}_{\epsilon}\} < r <\mu_D$ for some $\epsilon>0$. We divide both sides of Inequality \eqref{eq:OptCondP_r} by $\sigma_D^2/(\mu_D-r)$ which implies that for all $\epsilon>0$
\begin{equation}
    \label{eq:r_cCond}
    \frac{r-\Tilde{r}}{\mu_D-r} < 1+\epsilon \quad \mbox{or} \quad \frac{r-\Tilde{r}}{\mu_D-r} \leq 1
\end{equation}
for all $r\in(\max\{\Tilde{r},\Tilde{r}_\epsilon\},\mu_D)$.
Inequality \eqref{eq:r_cCond} cannot hold for all $r\in(\max\{\Tilde{r},\Tilde{r}_\epsilon\},\mu_D)$ unless $\Tilde{r}=\mu_D$
which contradicts the assumption. 
\item Let $\hat{\C}: [0,\mu_D) \to \R_+$ be defined by $\hat{\C}(r) := h\sigma_D^2/\left(2(\mu_D-r)\right)+p(\mu_D-r)$.
Let $\hat{r}$ denote the minimizer of $\hat{\C}$, i.e., $\hat{r}_p:= \mu_D- \frac{\sqrt{h\sigma^2/2}}{\sqrt{p}}$. Now we have
\begin{equation}
    \label{eq:Lim_rhat}
    \lim_{p \to \infty} \frac{\C(C_{\hat{r}_p})}{\hat{\C}(\hat{r}_p)}=
    \lim_{p \to \infty}\frac{h\E[\Jinf(C_{\hat{r}_p})]+p (\mu_D - \hat{r}_p)}{h\sigma_D^2/\left(2(\mu_D-\hat{r}_p)\right)+p(\mu_D-\hat{r}_p)} = \lim_{p \to \infty}
    \frac{h\E[\Jinf(C_{\hat{r}_p})]+\sqrt{ph\sigma_D^2/2}}{\sqrt{2ph\sigma_D^2}}
    =1,
\end{equation}
where the first two equalities use definitions and algebra and the final equality follows from part \ref{lemma:Em_c}.
Similarly, by part \ref{lemma:rStar_limit} we obtain
\begin{equation}
    \label{eq:Lim_rp}
    \lim_{p \to \infty} \frac{\C(C_{r_p})}{\hat{\C}(r_p)}=
    \lim_{p \to \infty}\frac{h\E[\Jinf(C_{r_p})]+p (\mu_D - r_p)}{h\sigma_D^2/\left(2(\mu_D-r_p)\right)+p(\mu_D-r_p)}=1.
\end{equation}
Optimality of $r_p$ together with \eqref{eq:Lim_rp} implies $\lim_{p \to \infty} \frac{\C(C_{\hat{r}_p})}{\hat{\C}(r_p)} \geq 1$.
Combining this with \eqref{eq:Lim_rhat}, we obtain
\begin{equation}
    \label{eq:Lim_rhat_rp}
    \lim_{p \to \infty} \frac{\C(C_{\hat{r}_p})}{\hat{\C}(r_p)} = \lim_{p \to \infty} \frac{\sqrt{2h\sigma_D^2}}{h\sigma_D^2/\left(2\sqrt{p}(\mu_D-r_p)\right)+\sqrt{p}(\mu_D-r_p)}\geq1.
\end{equation}
Observe that $h\sigma_D^2/\left(2\sqrt{p}(\mu_D-r_p)\right)+\sqrt{p}(\mu_D-r_p) \geq \sqrt{2h\sigma_D^2}$ and $\sqrt{2h\sigma_D^2} < \infty$. This result together with \eqref{eq:Lim_rhat_rp} yields the result: $\lim_{p \to \infty} \sqrt{p}(\mu_D-r_p) = \sqrt{h\sigma_D^2/2}$.
\Halmos
\end{enumerate}

% \section{Proof of Lemma \ref{lemma:Jicv}}
% \label{section:Jicv}
% Define the sequence $\{\{\lC_t(P_{\xi}):= \left(X_t-J_t(P_{\xi})+q_t(P_{\xi})\right)^+\}_{t\in \N_0}\}_{\xi\geq0} $, and let $\lC_{[a,b]}(P_{\xi}):=\sum_{t=a}^b \lC_t(P_{\xi})$.

% \begin{lemma}
%     \label{lemma:QLOrdering}
%     \citep[Lemma 9.][]{vanJaarsveld2024}
%     Let $0 \leq \xi_1 \leq \xi_2$.
%     For any $\xb_0$, and $t\geq 1$,
%     \begin{enumerate}[(a)]
%         \item $q_{[1,t]}(P_{\xi_1}) \leq_{\text{icv}} q_{[1,t]}(P_{\xi_2})$,
%         \item $\lC_{[1,t]}(P_{\xi_2}) \leq_{\text{icx}} \lC_{[1,t]}(P_{\xi_1})$.
%     \end{enumerate}
% \end{lemma}
% \proof{Proof of Lemma \ref{lemma:Jicv}.}
%     Under a PIL policy with parameter $\xi\geq0$, by Equation \eqref{eq:JDynamics}, $J_{t+1} (P_{\xi}), t \in \N_0$ can be written by
%     \[
%         J_{t+1} (P_{\xi}) = \left(J_t (P_{\xi}) +q_t (P_{\xi})-D_t \right)^+  = J_t (P_{\xi}) +q_t (P_{\xi})-D_t + \lC(P_{\xi}).
%     \]
%     Using the telescopic sum we have
%     \[
%         J_{t+1} (P_{\xi}) = J_0 + q_{[1,t]}(P_{\xi}) - D_{[1,t]} - \lC_{[1,t]}(P_{\xi}).
%     \]
%     Combining this relation, Lemma \ref{lemma:QLOrdering}, and Theorem 4.A.1 \cite{shaked2007} provide the result.
%     \Halmos
% \endproof

\section{Proof of Lemma \ref{lemma:qbounded}}
% Let $\tilde{D}_t$ represent the demand realization at $t \in \N_0$.
We drop $P_{\xi}$ in $q_t(P_{\xi})$, for simplicity of notation. Then we use iteration \eqref{eq:JDynamics} $L$ times to find
% For any $t \geq L$ and $\xi \geq 0$, let
% $
%     J_{t-1} = \left(((J_{t-L}+q_{t-L}-D_{t-L})^+ +\dots )^+ + q_{t-1}-D_{t-1} \right)^+
% $.
%By Equation \eqref{eq:JDynamics}
\begin{alignat}{2}
    \label{eq:qDist}
    & q_{t+1} = \xi - \E[J_{t+1} \mid \xb_{t-L+1}] = \xi - \E_{D_{t-L+1},\dots,D_{t}}[(((J_{t-L+1}+q_{t-L+1}-D_{t-L+1})^+ +\dots)^+ + q_{t}-D_{t})^+] \nonumber \\
   & \leq \xi - \E_{D_{t-L+1},\dots,D_{t-1}}[J_t] - q_{t} + \mu_D
    = \mu_D + \E_{D_{t-L},\dots,D_{t-1}}[J_t]
    - \E_{D_{t-L+1},\dots,D_{t-1}}[J_t] .
\end{alignat}
The first equality holds since for any $a\in\R$, $a^+ \geq a$. The final equality follows from $q_t=\xi-\E[J_t\mid \xb_{t-L}]$.
% Next observe the subadditivity of $(\cdot)^+$ operator provides that $\E_{D_{t-L},\dots,D_{t-1}}[J_t] - \E_{D_{t-L+1},\dots,D_{t-1}}[J_t] \leq \E_{\Tilde{D}}[(D_{t-L}-\Tilde{D})^+]$, where $D_{t-L}$ and $\Tilde{D}$ are i.i.d. random variables.
Next observe that $\E[\E_{D_{t-L},\dots,D_{t-1}}[J_t]] = \E[\E_{D_{t-L+1},\dots,D_{t-1}}[J_t]]$. 
%Additionally, by definition $\E[J_t] = \xi - \E[q_t]$.
This observation combined with \eqref{eq:qDist}, and Lemma 1 of \cite{vanJaarsveld2024} imply the results.
% Note that $\xb_0$ is a recurrent state, thus the result holds for $q_L(P_\xi)$.
\Halmos

\section{Proof of Lemma \ref{lemma:deltaEVJqAsymtote}}
\label{section:deltaEVJqAsymtoteProof}
Define $\overline{\kappa}:= \sup_{r \in [0,\mu_D)} \kappa$.
Note that $\kappa$ is strictly positive and finite for all $r \in [0,\mu_D)$. This fact, combined with Lemma \ref{lemma:momentsPlus} ensures that $0<\overline{\kappa}<\infty$. Furthermore, $\sup_{r \in [0,\mu_D)} \mu_+ \leq 2 \overline{\kappa}$, and $g_r(x) \leq \overline{\kappa}$ for all $x\geq0$ by Lemma \ref{lemma:UplusFunc}.
% Let $\DB:=\max_{t\geq L} D_t$. By [TODO] $\E \left[\DB \right] <\infty$.
% By Lemma \ref{lemma:qbounded}, .
% For $\epsilon_1,\epsilon_2>0$, let $y_{\epsilon_1} \in \{y\geq0: g_r(y_1) \leq \epsilon_1, \forall y_1\geq y\}$ and $\xi_{\epsilon_1,\epsilon_2} \in \{\xi\geq0: \P\{J_t(P_{\xi_1}) \leq y_{\epsilon_1} +r\}\leq \epsilon_2, \forall t \in \N_0, \xi_1 \geq \xi\}$. By Lemma \ref{lemma:Junbounded}, $\xi_{\epsilon_1,\epsilon_2}<\infty$ exists for all $\epsilon_1,\epsilon_2>0$.
For $\xi \geq0$, let the random variable $I_{\xi}$ be given by:
\begin{equation*}
    I_{\xi}:=\frac{h \mu_+}{2(\mu_D-r)}\int_{J_t\left(P_{\xi}\right)}^{J_t\left(P_{\xi}\right)+q_t(P_{\xi})-r} g_r(y) dy.
\end{equation*}
% Let $\epsilon:= h \overline{\kappa} \left(\overline{\kappa} \epsilon_2 + \epsilon_1 (1-\epsilon_2) \right)$,
Then
\[
    \E[I_\xi]\leq \frac{h\mu_+}{2 (\mu_D - r)} \E[q_t(P_\xi)-r]\sup_xg_r(x) \leq h \overline{\kappa}^2
\]

% \begin{alignat}{2}
%     \label{eq:EI}
%     &\E[I_{\xi_{\epsilon_1,\epsilon_2}}] = \nonumber\\
%     &\E\left[I_{\xi_{\epsilon_1,\epsilon_2}} \vert J_t(P_\xi) \leq y_{\epsilon_1} +r \right]\P\{J_t(P_\xi) \leq y_{\epsilon_1} +r\} + \E\left[I_{\xi_{\epsilon_1,\epsilon_2}} \vert J_t(P_\xi) > y_{\epsilon_1} + r \right]\P\{J_t(P_\xi) > y_{\epsilon_1} +r\} \nonumber\\
%     &\leq \frac{h\mu_+ \vert \E[q_t(P_\xi)]-r \vert}{2(\mu_D-r)}\left(\overline{\kappa} \epsilon_2 + \epsilon_1 (1-\epsilon_2) \right) \leq h \overline{\kappa} \left(\overline{\kappa} \epsilon_2 + \epsilon_1 (1-\epsilon_2) \right) = \epsilon.
% \end{alignat}
\noindent
The first inequality follows from applying the mean value theorem to $I_{\xi}$, and the second inequality from Lemma \ref{lemma:qbounded}and from $\mu_+\leq 2\overline\kappa$.
% Since $\epsilon_1$ and $\epsilon_2$ are unrestricted, for any $\epsilon>0$ we can find the corresponding values of $\epsilon_1$ and $\epsilon_2$. Let $\xi_{\epsilon}:=\xi_{\epsilon_1, \epsilon_2}$ for some $\epsilon_1$ and $\epsilon_2$ corresponding to $\epsilon>0.$
% Notice that Lemma \ref{lemma:Junbounded} implies that for all $\xi \geq \xi_{\epsilon}$ Inequality \eqref{eq:EI} holds.
% Additionally, by Lemma \ref{lemma:COPheavyTraffic} and Equation \eqref{eq:xiAsymptGrowth}, for each $\xi_{\epsilon}$ there exists $p_{\epsilon}\geq0$ such that $\xi(r_p) \geq \xi_{\epsilon}$ for all $p \geq p_{\epsilon}$, which completes the proof.
\Halmos

\section{Proof of Lemma \ref{lemma:ImprovementLemma}}
\label{section:ImprovementLemmaProof}
We prove the result by induction.
The case $t_2=t_1$ holds by the definition of $v_r$ (cf. Definition \ref{def: valuefunc}).
Next assume the result holds for $t_2 \geq t_1$. Then, by Definition \ref{def: valuefunc},
\begin{alignat*}{2}
    &\E[v_r(J_{t_2+1})\vert J_{t_1}]= \E[\E[c_{t_2+1}(C_r) + v_r (J_{t_2+2})- \C(C_r) \vert J_{t_2+1}]\vert J_{t_1}]\\
    &= \E[c_{t_2+1}(C_r) + v_r (J_{t_2+2})\vert J_{t_1}] - \C(C_r).
\end{alignat*}
Plugging this relation into the induction hypothesis we have,
\[
    \E[c_{[t_1,t_2]} (C_r) \vert J_{t_1}] = v_r(J_{t_1}) - \E[c_{t_2+1}(C_r) + v_r (J_{t_2+2})\vert J_{t_1}] + (t_2+2-t_1) \C(C_r),
\]
which gives the result by algebraic rearrangement.
\Halmos

\end{APPENDIX}
%
%   or 
%
% Acknowledgments here
\ACKNOWLEDGMENT{%
The last author is supported by a VENI research talent grant from the Dutch Science Foundation (NWO).
% Enter the text of acknowledgments here
}% Leave this (end of acknowledgment)

% References here (outcomment the appropriate case) 

\bibliographystyle{informs2014} % outcomment this and next line in Case 1

\setlength{\bibsep}{0pt}
\renewcommand*{\bibfont}{\footnotesize}
\bibliography{bib} % if more than one, comma separated

\end{document}